\documentclass{article}
\usepackage{amsmath}
\usepackage{amssymb}
\usepackage{amsfonts}
\usepackage{graphicx}

\begin{document}

\title{An exact solution for a non-autonomous delay differential equation}


\author{Kenta Ohira\\
Future Value Creation Research Center,\\ Graduate School of Informatics, Nagoya University, Japan
}

\maketitle

\begin{abstract}
  We derive an exact solution for a simple non-autonomous delay differential equation (DDE) over the entire real-time axis, representing it as a sum of Gaussian-shaped dynamics with distinct peak positions. This marks the first explicit solution for non-autonomous DDEs and is a rare example even among general DDEs. The constructed solution offers key physical insights and facilitates the analysis of system properties, such as the envelope profile of the dynamics.
\end{abstract}

\section{Introduction}
In the study of dynamics with time delays, quite complex behaviors have been discovered. Time delays most commonly arise due to finite transmission and production times. Time delay is an intrinsic feature in many control and interaction systems and has been studied across various fields, such as mathematics, biology, physics, engineering, and economics (e.g., \cite{heiden1979, bellman1963, cabrera_1, hayes1950, insperger, kcuhler, longtinmilton1989a, mackeyglass1977, miltonetal2009b, ohirayamane2000, smith2010, stepan1989, stepaninsperger, szydlowski2010}).
``Delay differential equations"  (DDEs) are the main mathematical approach and modeling tool for such systems. In general, delays induce instability in stable fixed points, leading to oscillatory and more complex dynamics. Additionally, as the delay lengthens, the complexity of the dynamics increases. For example, the Mackey-Glass equation \cite{mackeyglass1977} shows a progression from monotonic convergence to transient oscillations, sustained oscillations, and finally chaotic dynamics as the delay parameter of the feedback function increases. With models like this, the pathway to complex behavior in many delay systems is gradually becoming understood (e.g., \cite{taylor}), but due to its complexity, there is still much room for further investigation in this field. 

Recently, we have been studying delay differential equations, focusing specifically on a category known as non-autonomous DDEs, where the coefficients vary with time. 
These equations are generally considered more challenging to analyze than autonomous DDEs and are primarily investigated with a focus on stability using approximate or numerical methods \cite{Busenberg1984,Volz1986,Ming1990,Ford2002,Liu2006,Gyori2017,Lucas2018,Kuptsov2020,Herrera2022}.  In particular, we have been examining the following simple non-autonomous delay differential equation, where $a$, $b$ are real constants and the delay $\tau \geq 0$  \cite{kentaohira2022,kentaohira2023,kentaohira2024}:
\begin{equation}\label{DDE1} 
\frac{dX(t)}{dt}  + a t X(t) = b X(t-\tau) 
\end{equation}
This equation represents a slight extension of Hayes' equation, which has been extensively studied \cite{hayes1950}. Hayes' equation is given as follows, where $\alpha$ and $\beta$ are real constants. 
\begin{equation}
\frac{dX(t)}{dt} + \alpha X(t) = \beta X(t-\tau) \label{Hayes} 
\end{equation} 
The difference between equations (\ref{DDE1}) and (\ref{Hayes}) lies in the second term, which becomes a linear function of time in equation (1). While this modification may seem minor, it leads to significantly different behavior. Notably, we have observed that the dynamics of equation (\ref{DDE1}) exhibit a frequency resonance phenomenon\cite{kentaohira2022,kentaohira2023} .

In this paper, we focus on the non-autonomous DDE described by equation (\ref{DDE1}) and demonstrate that its exact solution over the entire real time axis can be explicitly expressed as a sum of Gaussian-shaped dynamics with different peak positions. To the best of the authors' knowledge, this is the first instance where such an explicit solution has been obtained for non-autonomous DDEs. Moreover, it is a rare result even within the domain of DDEs in general.

The constructed solution provides valuable physical insights and enables us to explain several properties of the system, such as the envelope profile of its dynamics.

\section{Solution Stability} 

We consider the DDE \eqref{DDE1} and its solution for $t \geq 0$ and $\tau \geq 0$. In this case, the solution is determined by the choice of the initial function $X(t) = \phi(t)$ for $t \in [-\tau, 0]$. Although solving it mathematically is difficult, its stability can be analyzed.

\paragraph{\boldmath{Case $a \neq 0$}}
As $t \to \infty$, the second term of the DDE dominates, so for large $t$, the DDE approximates:
\begin{equation}
\frac{dX(t)}{dt} = -a t X(t).
\end{equation}
Solving this equation gives:
\begin{equation}
X(t) = {\cal{C}} e^{- \frac{1}{2} a t^2},
\end{equation}
where ${\cal{C}}$ is an arbitrary constant.

Thus, for $a > 0$, $X(t)$ converges to $0$ as $t \to \infty$. Conversely, for $a < 0$, the solution diverges.

From the convergence behavior of the solution, we conclude that for $a > 0$, the solution $X(t)$ converges for any initial function $X(t) = \phi(t), t \in [-\tau, 0]$. In other words, for $a > 0$, $X(t)$ is globally asymptotically stable at $X = 0$. Conversely, for $a < 0$, stability is lost.

\paragraph{\boldmath{Case $a=0$ and $\tau>0$}}
In this case, the DDE reduces to:
\begin{equation}
\frac{dX(t)}{dt} = b X(t-\tau),
\label{eq:Hayes_special_case}
\end{equation}
which corresponds to equation \eqref{Hayes} with $\alpha = 0$. For $\tau > 0$, it is known that the general solution can be expressed using the Lambert $W$ function \cite{shinozaki,pusenjak2017} as follows:
\begin{equation}
X(t) = \sum_{k = -\infty}^{\infty} {\cal{C}}_{k} e^{\lambda_k t},
\quad \lambda_k = \frac{1}{\tau} W_k (b \tau),
\label{eq:solution_Hayes_special_case}
\end{equation}
where ${\cal{C}}_{k}$ is determined by the initial function $\phi(t)$ for $t \in [-\tau, 0]$.

Here, the $W$ function is a multi-valued complex function defined for $z \in \mathbb{Z}$ by:
\begin{equation}
z = W(z) e^{W(z)}.
\label{eq:Wfunction_definition}
\end{equation}
For integer $k$, the $k$-th branch of the $W$ function is denoted as $W_k$. It is known that this solution is globally asymptotically stable at $X = 0$ for $- \frac{\pi}{2 \tau} < b < 0$, whereas it is unstable otherwise.

\paragraph{\boldmath{Case $a=0$ and $\tau=0$}}
In this case, the DDE simplifies to:
\begin{equation}
\frac{dX(t)}{dt} = b X(t).
\label{eq:DDE1_a_0_tau_0}
\end{equation}
Solving this equation gives:
\begin{equation}
X(t)={\cal{C}}e^{bt}.
\label{eq:solution_a_0_tau_0}
\end{equation}
Thus, the system is globally asymptotically stable for $b < 0$ and unstable for $b \geq 0$.

When we take the limit $\tau \to 0^+$, the stability condition for solution \eqref{eq:solution_Hayes_special_case} yields the same result as that for  \eqref{eq:solution_a_0_tau_0}. 
Therefore, for the entire case $a=0$, the stability condition can be stated as $- \frac{\pi}{2 \tau} < b < 0$.

In summary, based on the above analysis, the stability of the solution $X(t)$ of the DDE \eqref{DDE1} at $X = 0$ for $t \geq 0$ and $\tau \geq 0$ is as follows. 
\begin{itemize}
\item $a > 0$: Globally asymptotically stable.
\item $a = 0$ and $b \in \left(- \frac{\pi}{2 \tau}, 0 \right)$: Globally asymptotically stable.
\item $a = 0$ and $b \notin \left(- \frac{\pi}{2 \tau}, 0 \right)$: Unstable.
\item $a < 0$: Unstable.
\end{itemize}

\section{Analysis Using the Fourier Transform}
In the following sections, we apply the Fourier transform to the differential equation (\ref{DDE1}) and subsequently perform an inverse transform\cite{kentaohira2024}. Based on the properties of the transform, the function obtained through this process coincides with the solution to the differential equation (\ref{DDE1}) over \( t \in \mathbb{R} \). 
In the following, we exclude the case of \( X(t) = 0 \) from considerations as trivial.

\subsection{Fourier Transformability}


\paragraph{\boldmath{Case $a<0$:}}
$X(t)$ is not Fourier transformable. In the previous section, we showed that $X(t)$ diverges as $t \to \infty$. 
Therefore, since $X(t)$ is not absolutely integrable over $t \in \mathbb{R}$, it is not transformable.

\paragraph{\boldmath{Case $a=0$ and $\tau\neq0$:}}
$X(t)$ is not Fourier transformable. For $t\geq0$, the solution was given by
\begin{equation}
    X(t) = \sum_{k = -\infty}^{\infty} \mathcal{C}_{k} e^{\lambda_k t},
    \quad \lambda_k = \frac{1}{\tau} W_k (b \tau)
\end{equation}
Since this form of the solution satisfies Equation (\ref{eq:Hayes_special_case}) for all $t \in \mathbb{R}$, the general solution for $t \in \mathbb{R}$ is given by this expression. Here, the coefficients $\mathcal{C}_{k}$ are determined such that $X(t)$ is a real function, i.e., $\operatorname{Im}(X(t))=0$. However, since this solution is not absolutely integrable, it is not transformable.

\paragraph{\boldmath{Case $a=0$ and $\tau=0$:}}
In this case, the function is not Fourier transformable. The DDE is given by
\begin{equation}
    \frac{dX(t)}{dt} = b X(t)
\end{equation}
which leads to the solution
\begin{equation}
    X(t) = \mathcal{C} e^{bt}
\end{equation}
Since this function is not absolutely integrable, it is not transformable.

\paragraph{\boldmath{Case $a>0$:}}
$X(t)$ is Fourier transformable. 
In the previous section, we showed that for $a>0$, $X(t)$ converges to $0$ as $t \to \infty$. By similar reasoning, $X(t)$ also converges to $0$ as $t \to -\infty$. 
Therefore, since $X(t)$ is absolutely integrable over $t \in \mathbb{R}$, it is transformable.

\paragraph{Conclusion:}
From the above, for $a>0$, the solution $X(t)$ of DDE (\ref{DDE1}) is Fourier transformable over $t \in \mathbb{R}$. Otherwise, it is not transformable.


\subsection{Derivation of the General Solution in Integral Form}
In the parameter range where Fourier transform is possible (\( a > 0 \)), we seek the general solution to the delay differential equation (\ref{DDE1}) over \( t \in \mathbb{R} \).
\vskip\baselineskip

\noindent
For \( \tau = 0 \), the differential equation becomes:
\begin{equation}
    \frac{dX(t)}{dt} = (b - a t) X(t)
\end{equation}
Solving this yields with the integration constant ${\cal{C}}$:
\begin{equation}
    X(t) = {\cal{C}} e^{- \frac{1}{2}a t^2 + b t}
\end{equation}
which represents a Gaussian centered at \( t = b/a \).
\vspace{1em}\\
For \( \tau > 0 \), applying the Fourier transform to equation (\ref{DDE1}) gives a differential equation for \( \hat{X}(\omega) \):
\begin{equation}
    i \omega \hat{X}(\omega) + i a \frac{d \hat{X}(\omega)}{d \omega} = b e^{-i \tau \omega} \hat{X}(\omega)
\end{equation}
Solving this equation provides:
\begin{equation}
    \hat{X}(\omega) = {\cal{C}} \exp \left[ -\frac{1}{2a} \omega^2 + \frac{b}{a \tau} e^{-i \tau \omega} \right]
\end{equation}
By taking the inverse Fourier transform of \( \hat{X}(\omega) \), we obtain the general solution \( X(t) \) for the original differential equation:
\begin{equation}
    X(t) = \frac{1}{2\pi} \int_{-\infty}^{\infty} \hat{X}(\omega) e^{i t \omega} d\omega = \frac{{\cal{C}}}{2\pi} \int_{-\infty}^{\infty} \exp \left[ -\frac{1}{2a} \omega^2 + \frac{b}{a \tau} e^{-i \tau \omega} + i t \omega \right] d\omega
\end{equation}
This function \( X(t) \) satisfies the original differential equation (\ref{DDE1}) for any \( t \in \mathbb{R} \).

\subsection{Complex Integration}
In order to investigate the above solution, we proceed with additional calculations by completing the square for the terms involving \( \omega^2 \) and \( \omega \):
\begin{equation}
    X(t) = \frac{{\cal{C}}}{2\pi} e^{-\frac{1}{2}a t^2} \int_{-\infty}^{\infty} \exp \left[ -\frac{1}{2} \left( \frac{\omega}{\sqrt{a}} - i \sqrt{a} t \right)^2 + \frac{b}{a \tau} e^{a \tau t - i \sqrt{a} \tau s} \right] d\omega
\end{equation}
Using a variable change \( s = \frac{\omega}{\sqrt{a}} - i \sqrt{a} t \), we rewrite:
\begin{equation}
    X(t) = {\frac{{\cal{C}}}{2\pi}} \sqrt{a} e^{-\frac{1}{2}a t^2} \int_{-\infty - i \sqrt{a} t}^{\infty - i \sqrt{a} t} \exp \left[ -\frac{1}{2} s^2 + \frac{b}{a \tau} e^{a \tau t - i \sqrt{a} \tau s} \right] ds
\label{DDE1complex_integral_t=0}
\end{equation}

Let \( g(s) \) denote the integrand, then 
\begin{equation}
X(t) = {\frac{{\cal{C}}}{2\pi}}\sqrt{a}  e^{-{\frac{1}{2}}at^2} \int_{-\infty-i\sqrt{a}t}^{\infty-i\sqrt{a}t}g(s)ds.
\label{DDE1complex_integral_ts}
\end{equation}

\begin{figure}[ht]
    \hspace{-3cm}
\begin{center}
    \includegraphics[height=4.8cm]{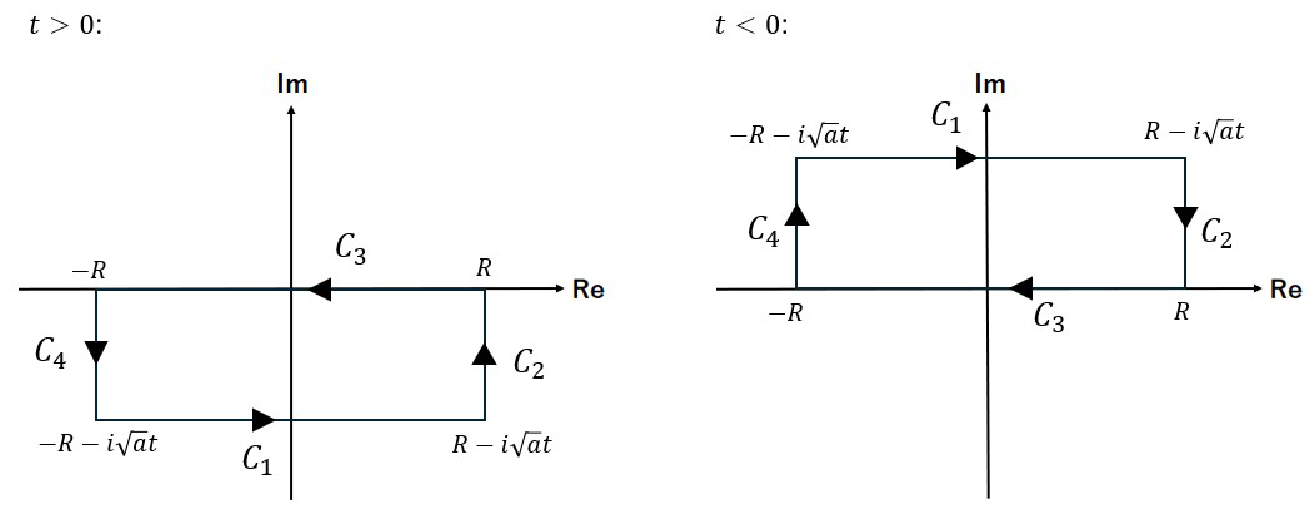}
    \caption{The integration path is defined as the rectangle in the complex plane with vertices at $z = -R - i\sqrt{a}t$, $R - i\sqrt{a}t$, $R$, and $-R$. The path on the left corresponds to $t > 0$, while the one on the right corresponds to $t < 0$. The direction of integration is indicated by the arrows.}
\end{center}
    \label{figure_DDE1_integral_path}
\end{figure}

For $t \neq 0$, we use the integration path shown in Figure \ref{figure_DDE1_integral_path} for complex integration. 
Since $g(s)$ is holomorphic over the entire complex plane, it satisfies:
\begin{equation}
    \lim_{R \to \infty} \oint_{C_1 + C_2 + C_3 + C_4} g(s) \, ds = 0.
\end{equation}
By considering the absolute value of the integral along the path $C_2$, the following inequality holds:
\begin{eqnarray}
    & \left| \int_{C_2} g(s) \, ds \right| = 
    \left| \int_{R - i\sqrt{a}t}^{R} \exp\left[-\frac{1}{2}s^2 + \frac{b}{a\tau} e^{a \tau t - i\sqrt{a}\tau s}\right] \, ds \right|\nonumber \\
    & \leq \sqrt{at} \exp\left[-\frac{1}{2}R^2 + \left|\frac{b}{a\tau}\right| e^{a\tau t}\right].
\end{eqnarray}
Since the right-hand side tends to $0$ as $R \to \infty$, we obtain:
\begin{equation}
    \lim_{R \to \infty} \int_{C_2} g(s) \, ds = 0.
\end{equation}
Similarly, performing the same calculation for the path $C_4$, we find:
\begin{equation}
    \lim_{R \to \infty} \int_{C_4} g(s) \, ds = 0.
\end{equation}
Therefore, the following holds:
\begin{equation}
    \lim_{R \to \infty} \int_{C_1} g(s) \, ds = -\lim_{R \to \infty} \int_{C_3} g(s) \, ds.
\end{equation}
This implies:
\begin{equation}
    \int_{-\infty - i\sqrt{a}t}^{\infty - i\sqrt{a}t} g(s) \, ds = \int_{-\infty}^{\infty} g(s) \, ds.
\end{equation}
This clearly also holds for the case of $t =0$.

Applying this result and Equation (\ref{DDE1complex_integral_ts}) to the expression for $X(t)$, we obtain:
\begin{equation}
    X(t) = \frac{{\cal{C}}}{2\pi} \sqrt{a} e^{-\frac{1}{2}at^2} 
    \int_{-\infty}^{\infty} \exp\left[-\frac{1}{2}s^2 + \frac{b}{a\tau} e^{a \tau t - i\sqrt{a}\tau s}\right] \, ds.
\end{equation}
The integral has been simplified, but solving it remains challenging.

\section{Series Solution}
Here, we guess the form of the solution based on the integral form obtained earlier and determine the solution by substituting it into the original DDE (\ref{DDE1}). Furthermore, we recognize that this solution remains well-defined even for $\tau<0$, leading to the general solution of the DDE for $a>0$, $\tau\neq0$, and $t \in \mathbb{R}$.

\subsection{Solution Hypothesis}
From the integral solution obtained in the previous section, we can hypothesize that
\begin{equation}
    X(t) = e^{-{1 \over 2} t^2} f(e^{a \tau t})
\end{equation}
where \( f : \mathbb{R} \rightarrow \mathbb{R} \) is a \( C^{\infty} \)-class function. Substituting this into the original differential equation
\begin{equation}
    \frac{dX(t)}{dt} + a t X(t) = b X(t - \tau),
\end{equation}
we find that \( f \) satisfies the following:
\begin{equation}
    f'(e^{a \tau t}) = \frac{b}{a \tau} e^{-{1 \over 2} a \tau^2} f(e^{a \tau (t - \tau)})
\end{equation}
Letting \( \alpha = e^{- a \tau^2} \) and \( \beta =   \frac{b}{a \tau} e^{-{1 \over 2} a \tau^2} \), \( f \) generally satisfies the following differential equation for \( x \in \mathbb{R} \):
\begin{equation}
    f'(x) = \beta f(\alpha x)
\end{equation}

\subsection{Derivation of the Taylor Series Expansion and Its Convergence}
While \( f(x) \) cannot be expressed using elementary functions, we can determine \( f^{(n)}(0) \) using it:
\begin{align}
    \begin{aligned}
        &f'(x) = \beta f(\alpha x), \quad f'(0) = \beta f(0) \\
        &f''(x) = \alpha \beta f'(\alpha x) = \alpha \beta^2 f(\alpha^2 x), \quad f''(0) = \alpha \beta^2 f(0) \\
        &f^{(3)}(x) = \alpha^2 \beta f''(\alpha x) = \alpha^3 \beta^3 f(\alpha^3 x), \quad f^{(3)}(0) = \alpha^3 \beta^3 f(0) \\
        &f^{(4)}(x) = \alpha^3 \beta f^{(3)}(\alpha x) = \alpha^6 \beta^4 f(\alpha^4 x), \quad f^{(4)}(0) = \alpha^6 \beta^4 f(0) \\
        &\dots \quad \dots \quad \dots \\
        &\dots \quad \dots \quad \dots \\
        &f^{(n)}(x) = \alpha^{\frac{n(n-1)}{2}} \beta^n f(\alpha^n x), \quad f^{(n)}(0) = \alpha^{\frac{n(n-1)}{2}} \beta^n f(0)
    \end{aligned}
\end{align}
Thus, setting \( {\cal{C}} = f(0) \), substituting values for \( \alpha \) and \( \beta \), and rearranging, \( f(x) \) can be expanded at \( x = 0 \) as follows:
\begin{equation}
    f(x) = {\cal{C}} \sum_{n = 0}^{\infty} \frac{1}{n!} \left( \frac{b}{a \tau} \right)^n \left( e^{-{1 \over 2} a \tau^2} \right)^{n^2} x^n
\end{equation}
The convergence of this series is demonstrated by the Cauchy-Hadamard theorem. Letting \( \gamma_n \) represent the coefficients of this series and \( r \) its radius of convergence:
\begin{eqnarray}
    & r = \frac{1}{\lim_{n \to \infty} \sqrt[n]{|\gamma_n|}} = \frac{1}{\lim_{n \to \infty} \sqrt[n]{\frac{1}{n!}} \left| \frac{b}{a \tau} \right| \left( e^{-{1 \over 2} a \tau^2} \right)^n} \nonumber \\
    & =
    \begin{cases} 
        0 & \text{if $a \leq 0$ or $\tau = 0$} \\
        \infty & \text{if $a > 0$  and $\tau \neq 0$}
    \end{cases}
\end{eqnarray}
Therefore, this series diverges for any \( x \in \mathbb{R} \) when \( a \leq 0 \) or \( \tau = 0 \), and converges for any \( x \in \mathbb{R} \) when \( a > 0 \) and \( \tau \neq 0 \). Thus, \( f \) is defined only for this range of parameters and delay. 

\subsection{Series Solution and Its Dynamics}
\begin{figure}[h]
\begin{center}
    \includegraphics[height=10cm]{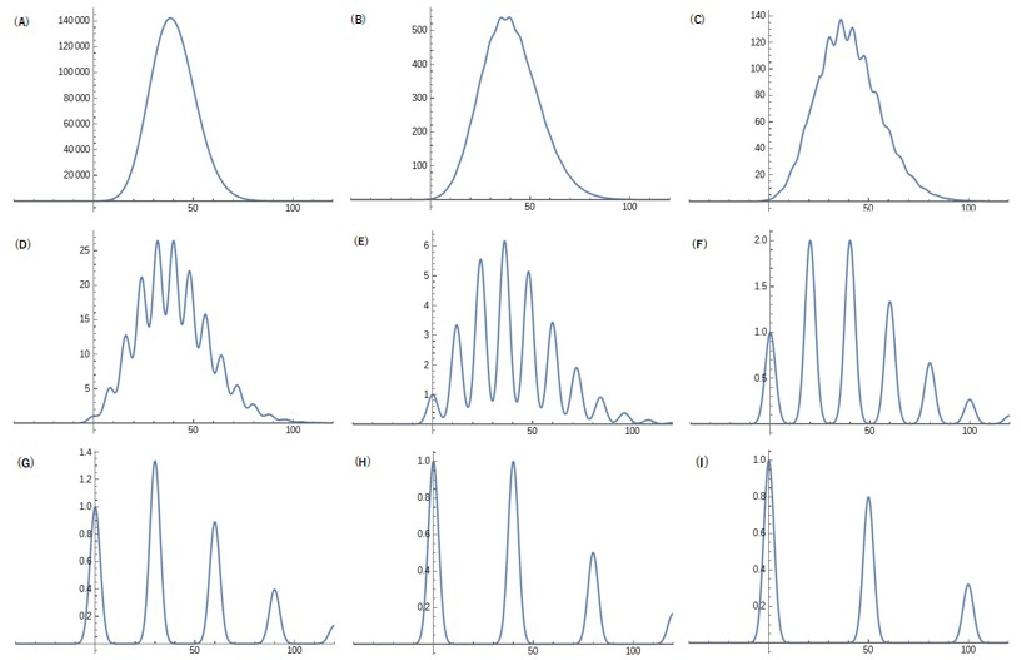}
    \caption{Graph of \( t - X(t) \) for the series solution (\ref{series_expansion_of_the_solution}) with \( a > 0 \) and \( \tau > 0 \), with varying delay \( \tau \). We set  \( {\cal{C}} =1 \) and parameters are \( a = 0.15 \), \( b = 6 \) with \( \tau \) set to (A) 3, (B) 5, (C) 6, (D) 8, (E) 12, (F) 20, (G) 30, (H) 40, and (I) 50. The series sums to \( n = 500 \).}
    \label{figure_DDE1_dynamics}
\end{center}
\end{figure}

\begin{figure}
\begin{center}
    \includegraphics[height=10cm]{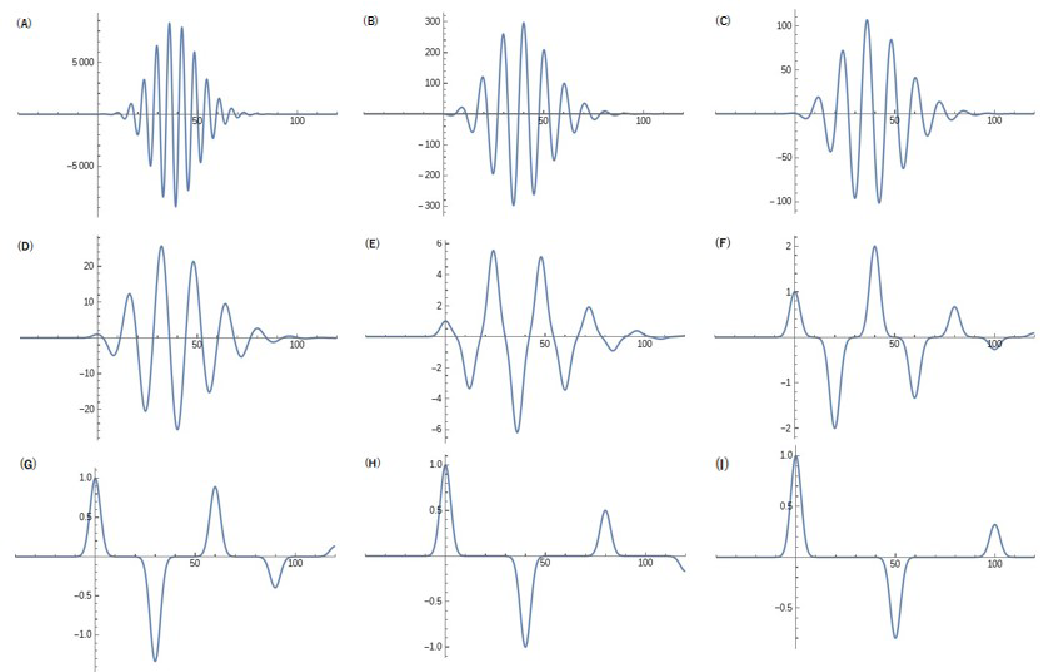}
    \caption{Graph of \( t - X(t) \) for the series solution (\ref{series_expansion_of_the_solution}) with \( a > 0 \) and \( \tau > 0 \), with varying delay \( \tau \).  We set  \( {\cal{C}} =1 \) and parameters are \( a = 0.15 \), \( b = -6 \) with \( \tau \) set to (A) 3, (B) 5, (C) 6, (D) 8, (E) 12, (F) 20, (G) 30, (H) 40, and (I) 50. The series sums to \( n = 500 \).}
    \label{figure_DDE1_dynamics_b_negative}
\end{center}
\end{figure}

From the above, when \( a > 0 \) and \( \tau > 0 \), the series solution \( X(t) \) to the DDE (\ref{DDE1}) is given by
\begin{equation}
    X(t) = e^{-{1 \over 2} a t^2} f(e^{a \tau t}) = {\cal{C}} e^{-{1 \over 2} a t^2} \sum_{n = 0}^{\infty} \frac{1}{n!} \left( \frac{b}{a \tau} \right)^n \left( e^{-{1 \over 2} a \tau^2} \right)^{n^2} \left( e^{a \tau t} \right)^n
    \label{series_expansion_of_the_solution0}
\end{equation}
where \( {\cal{C}} \) is an arbitrary constant such that \( f(0) = \lim_{t \to -\infty} f(e^{a \tau t}) = {\cal{C}} \). Further rearranging gives
\begin{equation}
    X(t) = {\cal{C}} \sum_{n = 0}^{\infty} \frac{1}{n!} \left( \frac{b}{a \tau} \right)^n e^{-{1 \over 2} a (t - n \tau)^2}
    \label{series_expansion_of_the_solution}
\end{equation}
This represents a superposition of Gaussians centered at \( t = n \tau \) with height \( \frac{{\cal{C}}}{n!} \left( \frac{b}{a \tau} \right)^n \). 
One can easily verify that  (\ref{series_expansion_of_the_solution}) is a solution of equation (\ref{DDE1}) through direct substitution.
Some representative example plots are shown in Figures \ref{figure_DDE1_dynamics} and \ref{figure_DDE1_dynamics_b_negative}.We observe that the series solution (\ref{series_expansion_of_the_solution}) aligns with the result previously obtained in integral form \cite{kentaohira2024} :
\begin{equation}
{X}(t) ={ {\cal{C}}\over{\sqrt {2 \pi a} }}\int_{-\infty}^{\infty} \exp \left[ {- {1\over 2 a} \omega^2 + {b\over \tau a} \cos (\omega \tau)}\right]\cos({b\over \tau a}\sin(\omega \tau) - \omega t )d\omega
\label{solx3}
\end{equation}
These results also  correspond to the dynamical trajectories obtained through direct numerical integration of Equation (\ref{DDE1}), using appropriate numerical initial functions as described in \cite{kentaohira2024}.

\vspace{3em}
\clearpage

\subsection{Extension to Negative Delay}
This solution is defined for \( a > 0 \) and \( \tau \neq 0 \), meaning it can be extended to both \( \tau > 0 \) and \( \tau < 0 \). \( C \) is an arbitrary constant.
\begin{equation}
X(t) = {\cal{C}} \sum_{n = 0}^{\infty} \frac{1}{n!} \left( \frac{b}{a \tau} \right)^n e^{-\frac{1}{2} a (t - n \tau)^2}, \quad {\cal{C}} = f(0) =
\begin{cases} 
\lim_{t \to -\infty} f(e^{a \tau t}) &\text{if \( \tau > 0 \)}\\
\lim_{t \to \infty} f(e^{a \tau t}) &\text{if \( \tau < 0 \)}
\end{cases}
\label{series_expansion_of_the_solution2}
\end{equation}
The difference between positive and negative \( \tau \) lies in the meaning of the arbitrary constant \( C \). For \( \tau > 0 \), setting \( C \) defines the "past" of \( f(e^{a \tau t}) \); for \( \tau < 0 \), it defines the ``future". Thus, \( \tau > 0 \) implies a ``delay," while \( \tau < 0 \) implies an ``advance."
\vskip\baselineskip

\noindent
Considering the relationship between the dynamics for \( \tau > 0 \) and \( \tau < 0 \), we set \( t = -t \), \( b = -b \), and \( \tau = -\tau \) in the original solution to obtain
\begin{equation}
X(t) = {\cal{C}} \sum_{n = 0}^{\infty} \frac{1}{n!} \left( \frac{b}{a \tau} \right)^n e^{-\frac{1}{2} a (-t + n \tau)^2}
\end{equation}
which is identical to the original solution. Therefore, if \( X(t) \) is considered a function of \( t \), \( b \), and \( \tau \) (denoted \( X(t, b, \tau) \)), then \( X(t, b, \tau) = X(-t, -b, -\tau) \), representing symmetry along the vertical \( X \)-axis. Some examples are shown in Figure 4.

\begin{figure}[h]
\hspace{-2cm}
\begin{center}
\includegraphics[height=6.0cm]{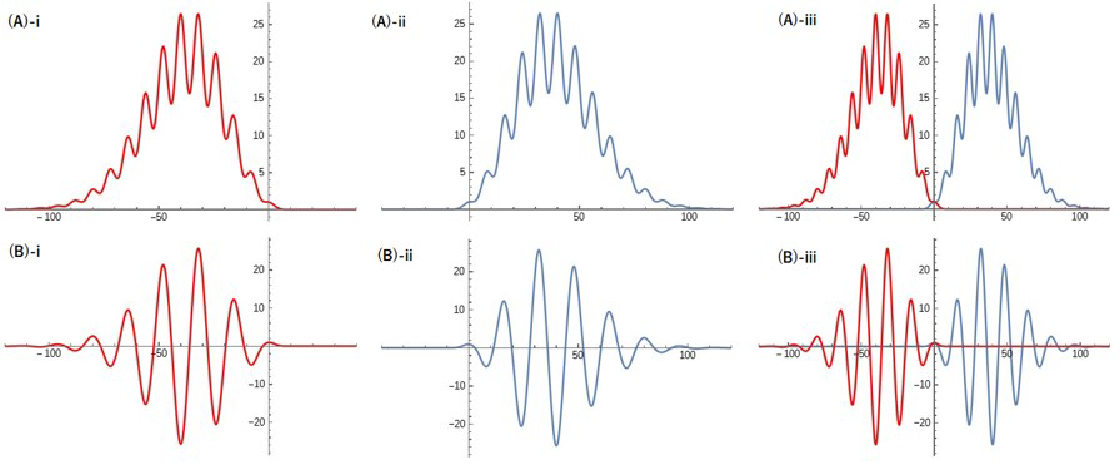}
\caption{For parameter values $(A)$ \( a = 0.15, b = -6 \) and $(B)$ \( a = 0.15, b = 6 \), with delay \( \tau = -8 \) in both cases. i: \( X(t, b, \tau) \), ii: \( X(t, -b, -\tau) \), iii: a comparison of both, showing symmetry.}
\end{center}
\label{figure_DDE1_dynamics_symmetry}
\end{figure}

\subsection{Solutions for Special Parameter Cases}
The solution of the DDE \eqref{DDE1} for \( a > 0 \) and \( \tau \neq 0 \), given in \eqref{series_expansion_of_the_solution0} or \eqref{series_expansion_of_the_solution}, is now used to study special parameter values. The solutions are derived as follows.

\paragraph{\boldmath{Case $a \to \infty$}}
In this case, the term with $n=0$ remains in the sum for the solution. Therefore,
\begin{equation}
    \lim_{a \to \infty}X(t) = \lim_{a \to \infty}{\cal{C}}e^{-\frac{1}{2} a t^2}.
\end{equation}
This leads to the following result:
\begin{equation}
    \lim_{a \to \infty}X(t) =
    \begin{cases}
        {\cal{C}}, & \text{if } t=0, \\
        0, & \text{if } t \neq 0.
    \end{cases}
    \label{eq:series_a_infty}
\end{equation}
In fact, in this case, the second term in the DDE \eqref{DDE1} becomes dominant:
\begin{equation}
    \frac{dX(t)}{dt} = -a t X(t).
\end{equation}
Thus, the solution follows as:
\begin{equation}
    X(t) = \lim_{a \to \infty}{\cal{C}}e^{-\frac{1}{2} a t^2},
\end{equation}
which matches the result above. 

\paragraph{\boldmath{Case $b=0$}}
In this case, only the term with $n=0$ in the sum remains, leading to the following solution:
\begin{equation}
    X(t) = {\cal{C}}e^{-\frac{1}{2} a t^2}.
    \label{eq:series_b_0}
\end{equation}
The DDE simplifies to:
\begin{equation}
    \frac{dX(t)}{dt} = -a t X(t),
\end{equation}
whose solution agrees with \eqref{eq:series_b_0}. 

\paragraph{\boldmath{Case $|b| \to \infty$}}
In this case, the form of the solution remains the same. 

\paragraph{\boldmath{Case $|\tau| \to \infty$}}
In this case, only the term with $n=0$ in the sum remains, leading to:
\begin{equation}
    \lim_{|\tau| \to \infty}X(t) = {\cal{C}}e^{-\frac{1}{2} a t^2}.
    \label{eq:series_tau_infty}
\end{equation}
Since \( \lim_{|t| \to \infty}X(t) = 0 \), the term \( bX(t-\tau) \) in the DDE vanishes. Hence,
\begin{equation}
    \frac{dX(t)}{dt} = -a t X(t).
\end{equation}
Solving this yields \eqref{eq:series_tau_infty}.

\subsection{Approximation of the Maximum}

Let \( X_{n}(t) \) be the Gaussian for each \( n \) in the series of (\ref{series_expansion_of_the_solution}), defined as
\begin{equation}
X_{n}(t) = \frac{{\cal{C}}}{n!} \left( \frac{b}{a \tau} \right)^n e^{-\frac{1}{2} a (t - n \tau)^2}
\end{equation}
Thus, \( X(t) \) can be written as
\begin{equation}
X(t) = \sum_{n = 0}^{\infty} X_{n}(t)
\end{equation}
In the following, we compare \( X(t) \) with the curve \( G(t) \), which encloses the independent arrangement of \( X_{n}(t) \) (excluding the \( t \)-axis). As shown in Figures 5 and 6, by increasing \( |\tau| \) while keeping \( a \) and \( b \) fixed, \( G(t) \) approaches \( X(t) \) since larger \( |\tau| \) reduces the mutual influence among each Gaussian \( X_n(t) \) in \( X(t) \).

The absolute value of each \( X_{n}(t) \) reaches its maximum value at \( t = n\tau \), given by \( \frac{{|\cal{C}|}}{n!} \left|\frac{b}{a \tau}\right|^n \), so the maximum value of \( |G(t)| \) is given by:
\begin{equation}
\underset{t \in \mathbb{R}}{max}(|G(t)|)=max\{|X_n(n\tau)|\}_{n=0,1,2,...}=max\left\{\frac{{|\cal{C}|}}{n!} \left|\frac{b}{a \tau}\right|^n\right\}_{n=0,1,2,...}
\end{equation}
Furthermore, if \( n = n^{*} \) gives the maximum value \( \frac{{|\cal{C}|}}{n!} \left|\frac{b}{a \tau}\right|^n \), then:
\begin{equation}
\underset{t \in \mathbb{R}}{max}(|G(t)|)=|G(n^{*}\tau)|= |X_{n^{*}}(n^{*}\tau)| = \frac{{|\cal{C}|}}{n^{*}!} \left|\frac{b}{a \tau}\right|^{n^{*}}
\end{equation}
Also, $n^{*}$ satisfies the following:
\begin{equation}
	n^{*}
	\begin{cases}
		= 0 \quad &\text{if } |{b\over a\tau}|<1 \\
		=0,1 \quad &\text{if } |{b\over a\tau}|=1 \\
		\geq 1 \quad &\text{if } |{b\over a\tau}|>1
	\end{cases}
\end{equation}
Since \( X(t) \) is a superposition of Gaussians with different centers, it can be assumed that the maximum value of \( |X(t)| \) occurs near \( t = n^{*} \tau \), where \( |G(t)| \) is also maximized:
\begin{equation}
\underset{t \in \mathbb{R}}{max}(|X(t)|)\approx |X(n^{*}\tau)|
\end{equation}

Regarding \( n^{*} \), when \( \left|\frac{b}{a \tau}\right| < 1 \), there is only one solution at \( n^{*} = 0 \). When \( \left|\frac{b}{a \tau}\right| \geq 1 \), \( n^{*} \geq 1 \) always exists. If two values of \( n^{*} \) exist, or if there exists an \( |X_{n}(n\tau)| \) very close to \( |X_{n^{*}}(n^{*}\tau)| \), this method might make it difficult to estimate the time $t$ that gives the maximum value of \( |X(t)| \).
\vspace{2em}

\begin{figure}[ht]
\begin{center}
\includegraphics[height=6.5cm]{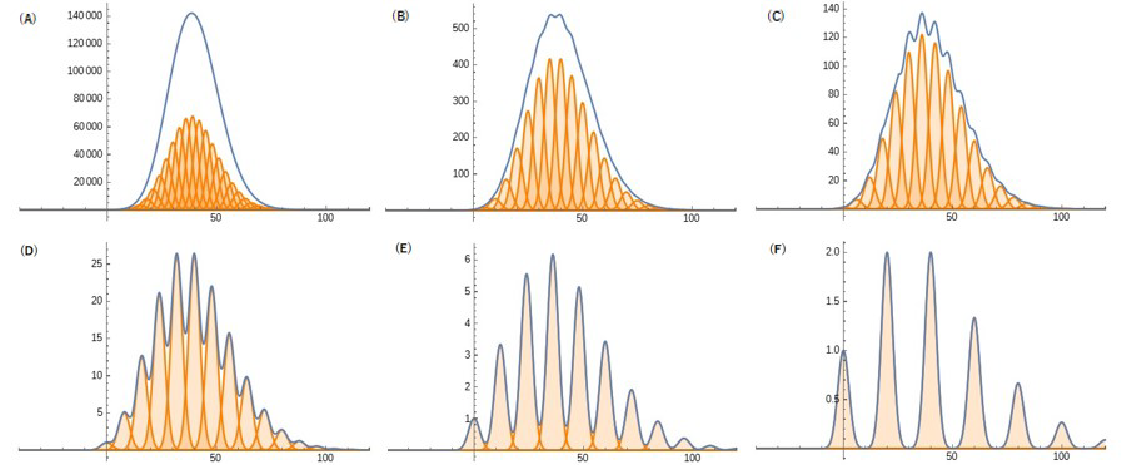}
\end{center}
\label{figure_DDE1_Dynamics_Gaussian}
\caption{
Dynamics of $X(t)$ (blue) and $G(t)$ (orange). $G(t)$ is defined as the curve surrounding the figure. We set  \( {\cal{C}} =1 \) and the parameters are \( a = 0.15, b = 6.0 \). \( \tau \) takes the values (A) 3, (B) 5, (C) 6, (D) 8, (E) 12, (F) 20.
}
\end{figure}

\clearpage

\begin{figure}
\begin{center}
\includegraphics[height=6.5cm]{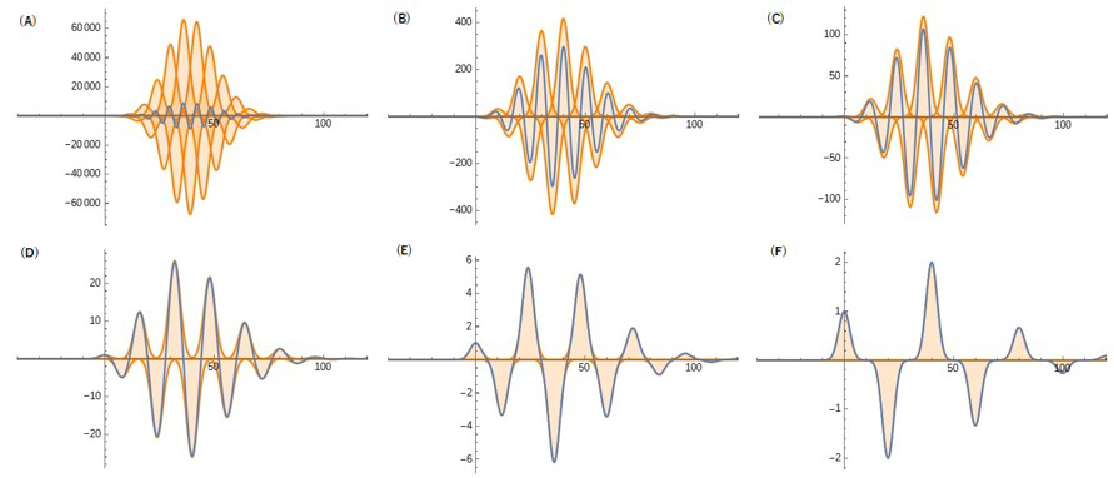}
\label{figure_DDE1_Dynamics_Gaussian_b_negative}
\end{center}
\caption{
Dynamics of $X(t)$ (blue) and $G(t)$ (orange). $G(t)$ is defined as the curve surrounding the figure.  We set  \( {\cal{C}} =1 \) and the parameters are  \( a = 0.15, b = -6 \). \( \tau \) takes the values (A) 3, (B) 5, (C) 6, (D) 8, (E) 12, (F) 20.
}
\end{figure}

\subsection{Approximation of the Envelope}

\( X_n(n\tau) = \frac{{\cal{C}}}{n!} \left(\frac{b}{a \tau}\right)^n \) is discrete with respect to \( n \). Extending it to a continuous function by letting \( n = t/\tau \) and defining it as \( E(t) \), then \( E(n\tau) = X_n(n\tau) \) forms the desired curve. \( E(t) \) is given by the following expression. (The case where \( b/\tau = 0 \) is excluded, as \( X(t) = {\cal{C}} e^{-\frac{1}{2} a t^2} \) results in a single point \( X_0(0) = {\cal{C}} \)):
\begin{equation}
E(t) =
\begin{cases} 
\frac{{\cal{C}}}{\Gamma\left(\frac{t}{\tau} + 1\right)} \left(\frac{b}{a \tau}\right)^{\frac{t}{\tau}} & \text{if \( b/\tau > 0 \)}\\
\pm \frac{{\cal{C}}}{\Gamma\left(\frac{t}{\tau} + 1\right)} \left(-\frac{b}{a \tau}\right)^{\frac{t}{\tau}} & \text{if \( b/\tau < 0 \)}
\end{cases}
\end{equation}
This behaves as an envelope for the waveform in the dynamics of \( G(t) \). Therefore, when \( \tau \) is sufficiently large, \( E(t) \) can approximate the waveform envelope of \( X(t) \). The limit of \( E(t) \) as \( b/\tau \) approaches zero is given by:
\begin{align}
\begin{aligned}
&\lim_{b/\tau \to +0} E(t) =
\begin{cases}
{\cal{C}} & \text{if \( t = 0 \)} \\
0 & \text{if \( t \neq 0 \)}
\end{cases}\\
&\lim_{b/\tau \to -0} E(t) =
\begin{cases}
\pm {\cal{C}} & \text{if \( t = 0 \)} \\
0 & \text{if \( t \neq 0 \)}
\end{cases}
\end{aligned}
\end{align}
Some examples are shown in Figures 7 and 8.

\begin{figure}[h]
\begin{center}
\includegraphics[height=6cm]{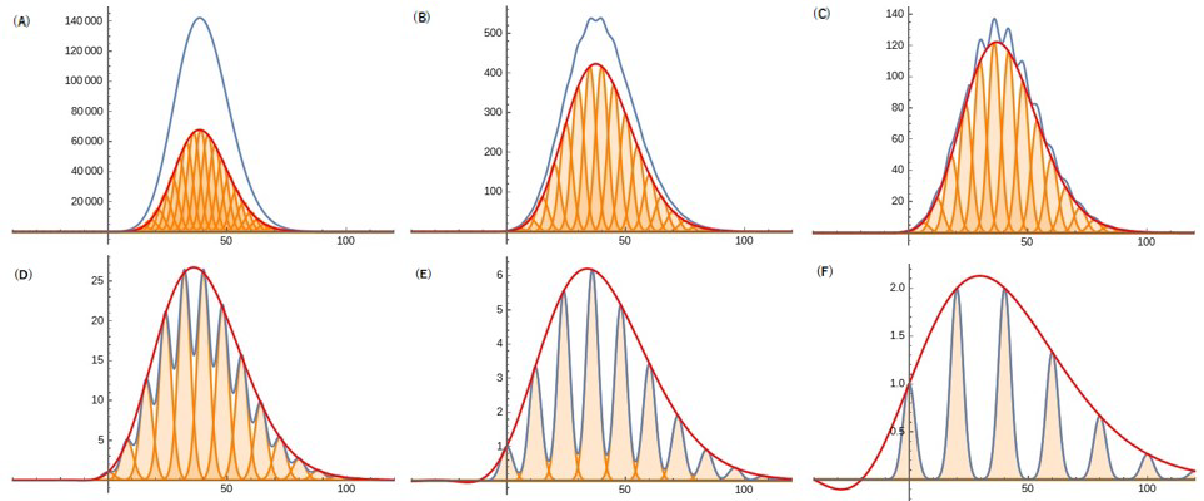}
\end{center}
\label{figure_DDE1_Dynamics_Gaussian_Envelope}
\caption{Dynamics of \( X(t) \) (blue), \( G(t) \) (orange), and \( E(t) \) (red). We set  \( {\cal{C}} =1 \). The parameters are \( a = 0.15, b = 6.0 \), and \( \tau \) takes the values (A) 3, (B) 5, (C) 6, (D) 8, (E) 12, (F) 20.}
\end{figure}

\begin{figure}
\begin{center}
\includegraphics[height=6cm]{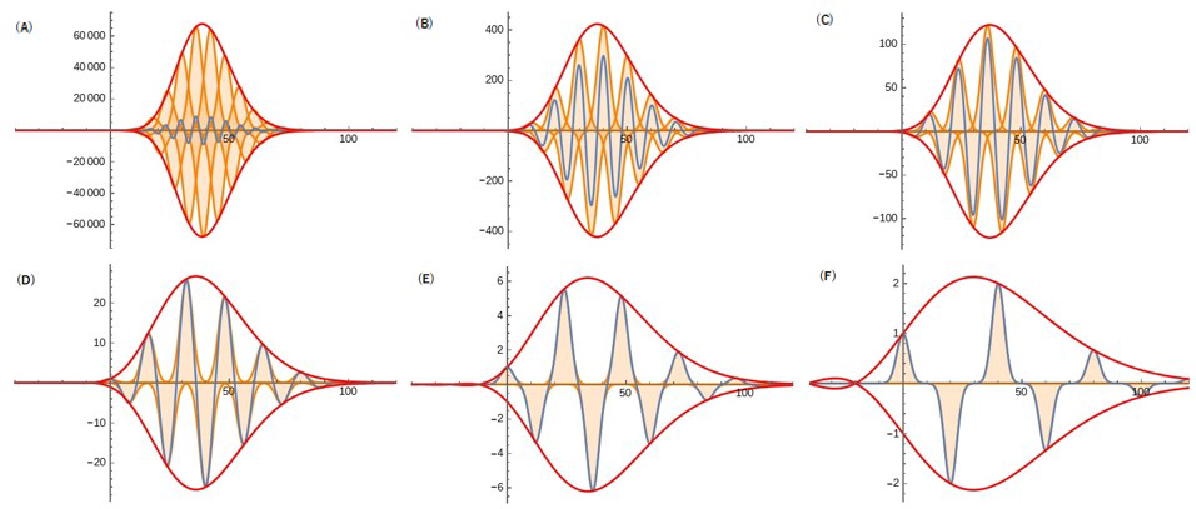}
\end{center}
\label{figure_DDE1_Dynamics_Gaussian_Envelope_b_negative}
\caption{Dynamics of \( X(t) \) (blue), \( G(t) \) (orange), and \( E(t) \) (red). We set  \( {\cal{C}} =1 \). The parameters are  \( a = 0.15, b = -6 \), and \( \tau \) takes the values (A) 3, (B) 5, (C) 6, (D) 8, (E) 12, (F) 20.}
\end{figure}

\clearpage


\section{Summary of General Solutions and Discussion}
To summarize, the solutions of the differential equation
\begin{equation}
\frac{dX(t)}{dt} + a t X(t) = b X(t-\tau),\quad t \in \mathbb{R},\quad a,b,\tau \in \mathbb{R}
\end{equation}
can be categorized as follows:
\vskip\baselineskip

\noindent
(i) For \( a = 0 \):
\begin{equation}
X(t) = \sum_{k = -\infty}^{\infty} C_k e^{\lambda_k t}, \quad \lambda_k = \frac{1}{\tau} W_k (b \tau)
\end{equation}
where \( C_k \) satisfies \( \operatorname{Im}(X(t)) = 0 \).
\vskip\baselineskip

\noindent
(ii) For \( \tau = 0 \):
\begin{equation}
X(t) = {\cal{C}} e^{- \frac{1}{2}a t^2 + b t}
\end{equation}
\vskip\baselineskip

\noindent
(iii) For \( a > 0 \) and \( \tau \neq 0 \):
\begin{equation}
X(t) = {\cal{C}} \sum_{n = 0}^{\infty} \frac{1}{n!} \left(\frac{b}{a\tau}\right)^n e^{-\frac{1}{2}a(t - n\tau)^2}
\end{equation}
\vskip\baselineskip

\noindent
(iv) For \( a < 0 \) and \( \tau \neq 0 \):
Currently unknown. However, as \( |t| \rightarrow \infty \), it diverges similarly to \( e^{-\frac{1}{2}a t^2} \).
\vskip\baselineskip

\noindent
Some discussion points are in order:

1. Our constructed solution is valid over the entire real-time axis, as mentioned. Given the uniqueness of the solution for first-order equations and the general unique correspondence of Fourier transform pairs, we believe that this solution is unique up to a scaling factor. However, a rigorous mathematical proof of this claim is still required.

In the typical formulation of DDEs, the initial function is specified for 
$t \in [-\tau, 0]$, and the solution is then sought for $t > 0$. When the initial function over 
$t \in [- \tau, 0]$ does not deviate significantly from the constructed solution, our result can serve as an approximation for such initial-value problems (see \cite{kentaohira2024} for details).

Overall, we hope that our exact solution for a non-autonomous DDE provides further insights into the nontrivial behaviors induced by delays.

\section*{Acknowledgments}
The author would like to thank Prof. Hideki Ohira and the members of his research group at Nagoya University for their valuable discussions. He also extends his gratitude to Prof. Prof. Yukihiko Nakata of Aoyama Gakuin University, Junya Nishiguchi of Tohoku University and Toru Ohira of Nagoya University for their constructive comments. This work was supported by the "Yocho-gaku" Project sponsored by Toyota Motor Corporation, the JSPS Topic-Setting Program to Advance Cutting-Edge Humanities and Social Sciences Research (Grant Number JPJS00122674991), JSPS KAKENHI (Grant Number 19H01201), and the Research Institute for Mathematical Sciences, an International Joint Usage/Research Center at Kyoto University.


\end{document}